\documentclass{elsart}
\usepackage{amssymb}
\usepackage{amsmath}
\usepackage[english]{babel}

\newtheorem{theorem}{Theorem}

\newtheorem{corollary}[theorem]{Corollary}
\newtheorem{definition}[theorem]{Definition}

\newtheorem{lemma}[theorem]{Lemma}

\newtheorem{remark}[theorem]{Remark}

\begin{document}
\begin{frontmatter}
\title{Fundamental solutions for a class of non-elliptic homogeneous differential
operators.}
\author{Brice Camus}
\address{Ruhr-Universit\"at Bochum, Fakult\"at f\"ur Mathematik,\newline%
Universit\"atsstr. 150, D-44780 Bochum, Germany.\newline E-mail:
brice.camus@uni-duisburg-essen.de}
\begin{abstract}
We compute temperate fundamental solutions of homogeneous
differential operators with real-principal type symbols. Via
analytic continuation of meromorphic distributions, fundamental
solutions for these non-elliptic operators can be constructed in
terms of radial averages and invariant distributions on the unit
sphere.
\end{abstract}
\end{frontmatter}
\section{Introduction and main results.}
If $P:=P(D_x)$ is a differential operator on $\mathbb{R}^n$ a
temperate fundamental solution to $P$ is a distribution
$\mathfrak{s}\in \mathcal{S}'(\mathbb{R}^n)$ such that $P(D_x)
\mathfrak{s} =\delta$, where $\delta$ is the delta-Dirac
distribution at the origin. Fundamental solutions play a major
role in the theory of PDE. For a large overview on this subject,
and applications, we refer to \cite{HOR1} vol. 1 \& 2. It is well
known, see e.g. \cite{Bjo,G-C,HOR1}, that differential operators
with constant coefficients have temperate fundamental solutions.
But, apart in very trivial cases like the Laplacian, it is
difficult to produce explicitly a solution. The case of order 3
homogeneous operators, in dimension 3, was treated in \cite{Wag}.
Always in dimension 3, the case of elliptic quartic operators was
considered in \cite{Wag1} and our contribution in \cite{Cam} was
to obtain temperate fundamental solutions for homogeneous elliptic
operators of any degree and in any dimension. Also, we mention
that the book of J.E. Bj\"ork \cite{Bjo} contains a very nice
study of the algebraic and analytic properties of fundamental
solutions for operators with polynomial or analytic symbols and
constant coefficients. In particular the presence of logarithmic
distributions, as occurring
in the present contribution, is predicted in a very general setting.\medskip\\
\textbf{Hypotheses and definitions.}\\
We are here interested in the case of a non-definite homogenous
polynomial $p$ on $\mathbb{R}^n$, i.e., $p(\lambda
\xi)=\lambda^{k} p(\xi)$. In all this article $k$ is the degree of
$p$. To simplify, we restrict our study to a real principal type
singularity, i.e. we assume that:
\begin{equation*}
(\mathcal{H})\,:\,\,\, \left\{
\begin{matrix}
p \text{ is real valued},\\
p(x)=0 \text{ and }\nabla p(x) =0 \Leftrightarrow x=0.
\end{matrix}
\right.
\end{equation*}
But $p$ complex valued is admissible, see section 2. In what
follows, we write:
\begin{equation*}
\mathfrak{C}(p)=\{\theta \in \mathbb{S}^{n-1}\text{ / }
p(\theta)=0\},
\end{equation*}
the trace of the characteristic set of $p$ on the unit-sphere. In
terms of polar coordinates, $(\mathcal{H})$ implies that the
restriction of $p$ to $\mathbb{S}^{n-1}$ satisfies:
\begin{equation*}
\nabla_\theta p(\theta)\neq 0\text{ near } \mathfrak{C}(p).
\end{equation*}
By a standard result of differential geometry, see e.g. \cite{G-C}
chapter 3, condition $(\mathcal{H})$ insures the existence of a
canonical (n-2)-dimensional measure $d\mathfrak{L}$ smooth on the
level sets $p(\theta)=\varepsilon$, for $\varepsilon>0$ small
enough. This measure, traditionally called Liouville or
Guelfand-Leray measure, satisfies the coarea formula:
\begin{equation*}
\int\limits_{\mathbb{S}^{n-1}} h(\theta) d\theta=
\int\limits_{\mathbb{R}}\left( \,\int\limits_{p(\theta)=u}
hd\mathfrak{L}\, \right) du,
\end{equation*}
for all $h$ with support in
$K_\varepsilon=\{\theta\in\mathbb{S}^{n-1}\text{ /
}|p(\theta)|\leq \varepsilon\}$. This relation defines a new
function: $u\mapsto \mathfrak{L}(h)(u)$, obtained by integration
of $f$ in the fibers $p^{-1}(u)$ w.r.t. $d\mathfrak{L}$. By Sard's
Theorem this function is finite almost everywhere and for any
$h\in C^\infty(\mathbb{S}^{n-1})$ it is easy to check that
$\mathfrak{L}(h)$ can be extended as an integrable function with
compact support
$\mathrm{supp}(\mathfrak{L}(h))\subset[\inf\limits_{\mathbb{S}^{n-1}}
p(\theta),\sup\limits_{\mathbb{S}^{n-1}} p(\theta)]$. With these
elementary facts in mind we introduce:
\begin{definition}
For a general function $g\in\mathcal{S}(\mathbb{R}^n)$ we define
the polar Guelfand-Leray transform of $g$ as:
\begin{equation*}
\mathfrak{L}(g(r\theta))(u):=\mathfrak{L}(g)(r,u)=\int\limits_{p(\theta)=u}
g(r\theta) d\mathfrak{L}(\theta),
\end{equation*}
simply by viewing the radius $r$ as a parameter.
\end{definition}
In all what follows the map $\mathfrak{L}$ is defined w.r.t. the
restriction of $p$ to $\mathbb{S}^{n-1}$ and:
\begin{equation*}
\mathfrak{L}^{(l)}(g(r\theta))(u)=\frac{d^l}{du^l}\left(\int\limits_{p(\theta)=u}
g(r\theta) d\mathfrak{L}(\theta)\right),
\end{equation*}
is the exterior derivative of degree $l$ w.r.t. the argument $u$.
Finally:
\begin{equation*}
\hat{h}(\xi)=\int\limits_{\mathbb{R}^n} e^{-i\langle x,\xi
\rangle} h(x)dx,
\end{equation*}
stands for the Fourier transform.\medskip\\
\textbf{Main results.}
\begin{theorem}\label{main}
Assume that $n\geq 2$ and that the symbol $p$ satisfies condition
$(\mathcal{H})$. A fundamental solution $\mathfrak{s}\in
\mathcal{S}'(\mathbb{R}^n)$ to $P$ is respectively given by:\medskip\\
\textbf{A)} If $k<n$ (locally integrable singularity):
\begin{equation*}
\left\langle \mathfrak{s}, f \right\rangle=
\frac{1}{(2\pi)^n}\int\limits_{0}^\infty \left\langle \log(|u|)\,
;\,\mathfrak{L}^{(1)} (\hat{f}(r\theta))(u) \right\rangle
r^{n-k-1} dr.
\end{equation*}
\textbf{B)} If $k\geq n$ (non-integrable case) then we have:
\begin{gather*}
\left\langle \mathfrak{s}, f \right\rangle=
\frac{1}{(2\pi)^n}\frac{\gamma+\Psi(k)}{\Gamma(2k)}
\frac{\partial^{2k-1}}{\partial r^{2k-1}} \left( r^{k+n-1}
\left\langle \log(|u|)\, ;\,
\mathfrak{L}^{(1)}(\hat{f}(r\theta))(u)\right\rangle
\right)_{|r=0}\\
+\frac{1}{(2\pi)^n\Gamma(1+2k)}\frac{\partial^{2k-1}}{\partial
r^{2k-1}} \left( r^{k+n-1} \left\langle \log(|u|)^2\, ;\,
\mathfrak{L}^{(1)}(\hat{f}(r\theta))(u)\right\rangle
\right)_{|r=0}\\
+\frac{1}{(2\pi)^n\Gamma(k)}\int\limits_{\mathbb{R}^{+}}
\log(r)\frac{\partial^{2k}}{\partial r^{2k}} \left( r^{k+n-1}
\left\langle \log(|u|)\, ;\, \mathfrak{L}^{(1)}(\hat{f}(r
\theta))(u)\right\rangle \right) dr.
\end{gather*}
Here $\gamma$ is Euler's constant and $\Psi(z)=\Gamma'(z)
/\Gamma(z)$.
\end{theorem}
The trivial case $n=1$, i.e. a monomial symbol, can be treated
directly and for $n=2$ the map $\mathfrak{L}$ is simply related to
Dirac masses at $\mathbb{S}^1\cap\{p=0\}$. Note that the results
are very different from the case of an elliptic operator. In
particular observe the presence of singularities supported in the
lacuna set of $p$ since distributions $\log(|u|)^j$, $j=1,2$ are
not smooth in $u=0$.

For non-integrable singularities we can say more and the method we
use allows to produce a one-parameter family of solutions:
\begin{corollary}\label{coro main}%
Under the conditions of Theorem \ref{main} and if $k\geq n$ a
temperate solution of $P(D)\mathfrak{s_0}=0$ is given by:
\begin{equation*}
\langle\mathfrak{s_0},f\rangle=\frac{\partial^{2k-1}}{\partial
r^{2k-1}} \left( r^{k+n-1} \left\langle \log(|u|)\, ;\,
\mathfrak{L}^{(1)}(\hat{f}(r\theta))(u)\right\rangle
\right)_{|r=0}.
\end{equation*}
Hence each $\mathfrak{s}+\lambda\mathfrak{s_0}$, $\lambda\in
\mathbb{C}$, is a temperate fundamental solution to $P(D)$.
\end{corollary}
\section{Proof of the main result.}
The strategy is as follows. If $p$ is positive for all $f\in
\mathcal{S}(\mathbb{R}^n)$ we have:
\begin{equation}\label{equ. of continuity}
\lim\limits_{\zeta\rightarrow 0} \frac{1}{(2\pi)^n}
\int\limits_{\mathbb{R}^n} p(\xi)^{\zeta}
\hat{f}(\xi)=f(0)=\left\langle \delta, f\right\rangle.
\end{equation}
If $p$ is a polynomial, or more generally an analytic function,
the integral in Eq.(\ref{equ. of continuity}) defines a
meromorphic distribution $\mathcal{P}(\zeta)$. See \cite{Bjo} for
this point. The Laurent development around $\zeta=-1$ can be
written:
\begin{equation}\label{laurent series}
\mathcal{P}(\zeta-1)=\sum\limits_{j=-1}^{-d}  \mu_j\zeta^{j}
+\mu_0 +\sum\limits_{j=1}^\infty \mu_j \zeta^j.
\end{equation}
But, according to Eq.(\ref{equ. of continuity}), we have:
\begin{equation*}
\lim\limits_{\zeta\rightarrow 0} \left\langle P(D)f,
\mathcal{P}(\zeta-1)\right\rangle =\left\langle \delta,
f\right\rangle,
\end{equation*}
and it follows that $\mu_0$ is a temperate fundamental solution to
$P(D)$.
\begin{remark}
\rm{Eq.(\ref{equ. of continuity}), combined with Eq.(\ref{laurent
series}), provides the set of relations:
\begin{equation*}
P(D)\mu_j =0,\text{ } \forall j <0,
\end{equation*}
in the sense of distributions of $\mathcal{S}'(\mathbb{R}^n)$. If
such non-zero terms exist, any affine combination
$\mu_0+\sum\limits_{j=1}^{d} \alpha_j\, \mu_{-j}$,
$(\alpha_1,...,\alpha_d)\in \mathbb{C}^d$, is a temperate
fundamental solution. This remark provides the basic strategy to
establish Corollary \ref{coro main}.}
\end{remark}
When $p$ is no more positive, or complex valued, the trick is to
compute the fundamental solution $\rho_0$ attached to $|p|^2$.
With $|p|^2=p(\xi) \bar{p}(\xi)$, it is easy to check that:
\begin{equation*}
\mu_0=\bar{P}(D) \rho_0,
\end{equation*}
is a fundamental solution to $P$. Hence, to attain our objective
we have to construct meromorphic extensions of the family of
distributions:
\begin{equation*}
\zeta\mapsto\int\limits_{\mathbb{R}^n} (|p(\xi)|^{2})^{\zeta}
g(\xi)d\xi,\, g\in \mathcal{S}(\mathbb{R}^n).
\end{equation*}
To solve a non-elliptic equation we transform the problem into a
positive, and hence simpler, problem. The expense is that
$|p(\xi)|^{-2}$ is more singular than $|p(\xi)|^{-1}$ and this
induces extra computations in the proof. We start by solving,
locally, the singularities of $p$. We have:
\begin{lemma} \label{normal forms}
If $p$ satisfies $(\mathcal{H})$ there exists local coordinates
$\omega$ (strictly speaking outside of the origin), such that we
have the local diffeomorphism:
\begin{equation*}
p(\xi)\simeq \left\{
\begin{matrix}
-\omega_1^k \text{ or } \omega_1^{k}, \text{ outside of } \mathfrak{C}(p)\times ]0,\infty[,\\
\omega_1^k \omega_2 \text{ in a neighborhood of }
\mathfrak{C}(p)\times ]0,\infty[.
\end{matrix}
\right.
\end{equation*}
\end{lemma}
\textit{Proof.} To blow up the singularity, we use polar
coordinates $\xi=(r,\theta)$. By homogeneity we have
$p(r\theta)=r^k p(\theta)$. First if $\theta_0 \notin
\mathfrak{C}(p)$ we choose:
\begin{equation}
(\omega_1,\omega_2,...,\omega_n)(r,\theta)=
(r|p(\theta)|^{\frac{1}{k}},\theta).
\end{equation}
We have $p(\xi)\simeq \pm \omega_1(r,\theta)^k$ in a conical
neighborhood of $\theta_0$. The sign is obviously given by the
sign of $p(\theta_0)$ and the Jacobian is $|J\omega|
(r,\theta)=|p(\theta)|^{\frac{1}{k}}\neq 0$. Next, if $\theta_0\in
\mathfrak{C}(p)$ by condition $(\mathcal{H})$ and by homogeneity
we have $\nabla_\theta p(\theta_0)\neq 0$. We can assume that
$\partial_{\theta_1} p(\theta) \neq 0$ and we chose:
\begin{equation*}
(\omega_1,\omega_2,\omega_3,...,\omega_n)(r,\theta)=
(r,p(\theta),\theta_2,...,\theta_{n_1}).
\end{equation*}
We have:
\begin{equation*}
|J\omega| (r,\theta_0)=|\frac{\partial p}{\partial
\theta_1}(\theta_0)|drd\theta \neq 0.
\end{equation*}
By continuity, this result holds in a sufficiently small
neighborhood of $\theta_0$. Since $\mathfrak{C}(p)$ is a compact
subset of $\mathbb{S}^{n-1}$ we can easily globalize the
construction.$\hfill{\blacksquare}$\medskip\\
To use these normal forms, we construct an adapted partition of
unity on $\mathbb{S}^{n-1}$. We pick a family of positive function
$\Omega_j$ on $\mathbb{S}^{n-1}$ such that:
\begin{equation*}
\sum\limits_{j=1}^N \Omega_j (\theta)=1 \text{ near }
\mathfrak{C}(p),
\end{equation*}
with the existence of a normal form $\omega_1^k \omega_2$ inside
each $\mathrm{supp}(\Omega_j)$. Next, since the previous
construction depends only on the set $\mathfrak{C}(p)$, we can
assume that $\mathrm{supp}(\Omega_j)\subset K_\varepsilon$ for
$\varepsilon>0$ chosen small enough so that the measures
$d\mathfrak{L}$ are well defined on each
$\mathrm{supp}(\Omega_j)$. Finally we can complete this finite set
as partition of unity on $\mathbb{S}^{n-1}$ with
$\Omega_0=1-\sum\limits_{j} \Omega_j$. The support of $\Omega_0$
is generally not connected, as shows the case $n=3$. With this
partition of unity we have:
\begin{equation*}
\int\limits_{\mathbb{R}^n} (|p(\xi)|^{2})^{\zeta} g(\xi)d\xi
=\sum\limits_{j=0}^N \,
\int\limits_{\mathbb{S}^{n-1}}\int\limits_{ \mathbb{R}^{+}}
\Omega_j(\theta) |p(r,\theta)|^{-2\zeta}
g(r,\theta)r^{n-1}drd\theta.
\end{equation*}
With this localization argument we use Lemma \ref{normal forms} to
trivialize locally the problem and we have to study the elementary
quantities:
\begin{gather*}
\mu^{\mathrm{ell}}(\zeta)=\int\limits_{\mathbb{R}^{+}} \omega_1^{2k\zeta} G(\omega_1) d\omega_1,\\
\mu_j^{\mathrm{sing}}(\zeta)=\int\limits_{\mathbb{R}^{+}\times\mathbb{R}}
\omega_1^{2k\zeta} (\omega_2^2)^{\zeta} G_j(\omega_1,\omega_2)
d\omega_1d\omega_2.
\end{gather*}
These new functions are obtained by pullback and integration:
\begin{gather*}
G(\omega_1)= \int \omega^{*} (\Omega_0(\theta)
g(r,\theta)r^{n-1})(\omega_1,...,\omega_n)d\omega_2...d\omega_n,\\
G_j(\omega_1,\omega_2)= \int \omega^{*} (\Omega_j(\theta)
g(r,\theta)r^{n-1})(\omega_1,...,\omega_n)d\omega_3...d\omega_n,
\end{gather*}
where $\omega^*$ stands for the pullback including the
multiplication by the Jacobian.\medskip\\
\textbf{Trivial contribution.}\\
We start by the analytic continuation of the elliptic part
$\mu^{\mathrm{ell}}(\zeta)$. We have:
\begin{equation*}
\frac{\partial^{2k}}{\partial \omega_1^{2k}}
\omega_1^{2k\zeta}=\omega_1^{2k(\zeta-1)}\prod\limits_{j=0}^{2k-1}
(2k\zeta-j) ,
\end{equation*}
and after $2k$ integrations by parts we obtain:
\begin{equation*}
\mu^{\mathrm{ell}}(\zeta-1)=\int\limits_{\mathbb{R}^{+}}
\omega_1^{2k(\zeta-1)} G(\omega_1)
d\omega_1=\left(\prod\limits_{j=0}^{2k-1}\frac{1}{(2k\zeta-j)}\right)
\int\limits_{\mathbb{R}^{+}} \omega_1^{2k\zeta}
\partial_{\omega_1}^{2k} G(\omega_1) d\omega_1.
\end{equation*}
The integral in the r.h.s. defines an holomorphic function near
$\zeta=0$. The constant term of the Laurent series at the origin,
determined by the rational function, is given by:
\begin{equation*}
\mu_0^{\mathrm{ell}}=\lim\limits_{\zeta\rightarrow 0}
\frac{\partial}{\partial\zeta} (\zeta\,
\mu^{\mathrm{ell}}(\zeta-1)).
\end{equation*}
With the holomorphic function near $\zeta=0$:
\begin{equation*}
h(\zeta)=\zeta\prod\limits_{j=0}^{2k-1}\frac{1}{
(2k\zeta-j)}=\frac{1}{2k} \prod\limits_{j=1}^{2k-1}
\frac{1}{(2k\zeta-j)},
\end{equation*}
we obtain:
\begin{equation*}
\mu_0^{\mathrm{ell}}=h'(0)\int\limits_{\mathbb{R}^{+}}
\partial_{\omega_1}^{2k}G(\omega_1)d\omega_1
+2k h(0)\int\limits_{\mathbb{R}^{+}} \log (\omega_1)
\partial_{\omega_1}^{2k} G(\omega_1) d\omega_1.
\end{equation*}
Clearly $2kh(0)=-1/\Gamma(2k)$ and a direct computation yields:
\begin{equation*}
h'(0)=-\frac{\gamma +\Psi(2k)}{\Gamma(2k)}.
\end{equation*}
Here $\Psi(\zeta)=\Gamma'(\zeta)/\Gamma(\zeta)$ is the usual
polygamma function of order 0 and
\begin{equation*}
\gamma=\lim\limits_{L\rightarrow \infty} \left(
\sum\limits_{j=1}^{L} \frac{1}{j}- \log(L)\right),
\end{equation*}
is Euler's constant.\medskip\\
\textbf{Non-trivial contribution.}\\
Now, we study the singular term
$\mu^{\mathrm{sing}}(\zeta)=\sum\limits_{j=1}^{N}\mu_j^{\mathrm{sing}}(\zeta)$.
We have:
\begin{gather*}
\frac{\partial^{2k+2}}{\partial\omega_1 ^{2k}\partial \omega_2^2}
\omega_1^{2k\zeta}
(\omega_2^2)^{\zeta}=\mathfrak{b}(\zeta)\omega_1^{2k(\zeta-1)}
(\omega_2^2)^{\zeta-1},\\
\mathfrak{b}(\zeta)=2\zeta(2\zeta-1)\prod\limits_{j=0}^{2k-1}
(2k\zeta-j) .
\end{gather*}
Accordingly, $\zeta=0$ is a pole of order 2 of the meromorphic
extension:
\begin{equation*}
\mu_j^{\mathrm{sing}}(\zeta-1)=\frac{1}{\mathfrak{b}(\zeta)}\int\limits_{\mathbb{R}^{+}\times\mathbb{R}}
\omega_1^{2k(\zeta-1)} (\omega_2^2)^{\zeta-1}\frac{\partial^{2k+2}
G_j}{\partial \omega_1^{2k}
\partial \omega_2^2} (\omega_1,\omega_2) d\omega_1 d\omega_2.
\end{equation*}
The constant term of the Laurent expansion is given by:
\begin{equation*}
\mu_{0,j}^{\mathrm{sing}}=\frac{1}{2}
\lim\limits_{\zeta\rightarrow 0} \frac{\partial^2}{\partial
\zeta^2} (\zeta^2 \,\mu_j^{\mathrm{sing}}(\zeta-1)).
\end{equation*}
Hence with the auxiliary functions:
\begin{gather*}
m(\zeta)=\frac{\zeta^2}{\mathfrak{b}(\zeta)}=\frac{1}{4k(2\zeta-1)}
\prod\limits_{j=1}^{2k-1} \frac{1}{2k\zeta-j},\\
M_j(\zeta)=\int\limits_{\mathbb{R}^{+}\times\mathbb{R}}
\omega_1^{2k(\zeta-1)} (\omega_2^2)^{\zeta-1}\frac{\partial^{2k+2}
G_j}{\partial \omega_1^{2k}
\partial \omega_2^2} (\omega_1,\omega_2) d\omega_1 d\omega_2,
\end{gather*}
we obtain that the term of interest is given by:
\begin{equation}
\mu_{0,j}^{\mathrm{sing}}=\frac{1}{2}
(m(0)M_j''(0)+2m'(0)M_j'(0)+m''(0)M_j(0)).
\end{equation}
By some elementary calculations we obtain respectively:
\begin{gather*}
m(0)=\frac{1}{2\Gamma(1+2k)},\\
m'(0)=\frac{1+k(\gamma +\Psi(2k))}{\Gamma(1+2k)}.
\end{gather*}
The coefficient $m''(0)$ plays no r\^ole here, see Eq.(\ref{coeff
zero}) below. The next step is to evaluate
$\mu_{0,j}^{\mathrm{sing}}$ in the coordinates $\omega$. After
integration by parts w.r.t. $\omega_2$, we have:
\begin{equation}\label{coeff zero}
M_j(0)=\int\limits_{\mathbb{R}^{+}\times\mathbb{R}}
\frac{\partial^{2k+2} G_j}{\partial \omega_1^{2k}
\partial \omega_2^2} (\omega_1,\omega_2) d\omega_1 d\omega_2=0.
\end{equation}
For the next distributional coefficient we find that:
\begin{gather}
M_j'(0)=\int\limits_{\mathbb{R}^{+}\times\mathbb{R}}(2k\log(\omega_1)+2\log(|\omega_2|))
\frac{\partial^{2k+2} G_j}{\partial \omega_1^{2k}
\partial \omega_2^2} (\omega_1,\omega_2) d\omega_1 d\omega_2 \notag \\
=\int\limits_{\mathbb{R}}2\log(|\omega_2|)) \frac{\partial^{2k+1}
G_j}{\partial \omega_1^{2k-1}
\partial \omega_2^2} (0,\omega_2) d\omega_2. \label{value M'}
\end{gather}
Finally, we obtain similarly:
\begin{gather}
M_j''(0)=\int\limits_{\mathbb{R}^{+}\times\mathbb{R}}(2k\log(\omega_1)+2\log(|\omega_2|))^2
\frac{\partial^{2k+2} G_j}{\partial \omega_1^{2k}
\partial \omega_2^2} (\omega_1,\omega_2) d\omega_1 d\omega_2 \notag\\
=4 \int\limits_{\mathbb{R}} \log(|\omega_2|)^2
\frac{\partial^{2k+1} G_j}{\partial \omega_1^{2k-1}
\partial \omega_2^2} (0,\omega_2) d\omega_2 \notag\\
+8k\int\limits_{\mathbb{R}^{+}\times\mathbb{R}}\log(\omega_1)\log(|\omega_2|)
\frac{\partial^{2k+2} G_j}{\partial \omega_1^{2k}
\partial \omega_2^2} (\omega_1,\omega_2) d\omega_1 d\omega_2. \label{value M''}
\end{gather}
After expanding the square in the integral we have, once more,
discarded the term attached to $\log(\omega_1)^2$, vanishing after
integration w.r.t. $\omega_2$.\medskip\\
\textbf{Invariant formulation.}\\
To achieve the proof we must formulate our distributions in a
geometrical way, also independent of the partition of unity
attached to the coordinates $\omega$. First, by construction, we
have to evaluate our distribution on $P(D)f$ so that after Fourier
transformation $g(\xi)=p(\xi)\hat{f}(\xi)$. Since $p$ is of degree
$k$, we have $G(\omega_1)=\mathcal{O}(\omega_1^{k+n-1})$ near
$\omega_1=0$. Same remark for
$G_{j}(\omega_1,\omega_2)=\mathcal{O}(\omega_1^{k+n-1})$ near
$\omega_1=0$. These properties are important since several
coefficient expressed below are related to Dirac-delta
distributions supported in $\omega_1=0$. According to
Eq.(\ref{value M'}) and Eq.(\ref{value M''}) at worst 3 different
terms occur which we treat separately distinguishing out the case
of $|p|^{-1}$ locally integrable or
not.\medskip\\
\textbf{1-Contribution of the elliptic directions.}\\
We have:
\begin{equation*}
\int\limits_{0}^{\infty}\partial_{\omega_1}^{2k}G(\omega_1)d\omega_1=-\partial_{\omega_1}^{2k-1}G(0).
\end{equation*}
If $2k-1<k+n-1$ this term vanishes and for $k\geq n$ we have:
\begin{equation*}
\partial_{\omega_1}^{2k-1}G(0)=
\frac{\partial^{2k-1}}{\partial r^{2k-1}}
\left(r^{n+k-1}\int\limits_{\mathbb{S}^{n-1}} \hat{f}(r\theta)
\Omega_0(\theta)\frac{d\theta}{p(\theta)}\right)_{|r=0}.
\end{equation*}
This identity holds after inversion of our diffeomorphism and the
substitution $g(\xi)=p(\xi)\hat{f}(\xi)$. When $k<n$, we can
integrate by parts the logarithmic contribution to obtain:
\begin{gather*}
\int\limits_{\mathbb{R}^{+}} \log (\omega_1)
\partial_{\omega_1}^{2k} G(\omega_1) d\omega_1=(2k-1)!\int\limits_{\mathbb{R}^{+}}
G(\omega_1) \frac{d\omega_1}{\omega_1^{2k}}\\
=(2k-1)!\int\limits_{\mathbb{R}^{+}\times\mathbb{S}^{n-1}}
\Omega_0(\theta) \hat{f}(r\theta) r^{n-k-1}
 dr\frac{d\theta}{p(\theta)}.
\end{gather*}
Observe that the integral w.r.t. $r$ is precisely convergent, for
any $f\in\mathcal{S}(\mathbb{R}^n)$, if and only if $k<n$. If
$k\geq n$ this argument does not holds, but we can write:
\begin{equation*}
\partial_{\omega_1}^{2k} G(\omega_1)=\frac{1}{2\pi}
\int\limits_{\mathbb{R}} e^{ip\omega_1} (ip)^{2k} \hat{G}(p)dp.
\end{equation*}
After inversion of our diffeomorphism and scaling out the
spherical term $p(\theta)$ in the phase, we obtain the
contribution:
\begin{equation*}
\int\limits_{\mathbb{R}^{+}} \log (\omega_1)
\partial_{\omega_1}^{2k} G(\omega_1) d\omega_1=%
\int\limits_{\mathbb{R}^{+}\times\mathbb{S}^{n-1}}\log(r|p(\theta)|^{\frac{1}{k}})
\frac{\partial^{2k}} {\partial r^{2k}}\left( \hat{f}(r\theta)
r^{n+k-1}\right)\Omega_0(\theta) dr\frac{d\theta}{p(\theta)}.
\end{equation*}
\textbf{2-Contribution of the non-elliptic directions.}\\
To express our amplitudes, we use the Schwartz kernel technique.
Let $\alpha=(\alpha_1,\alpha_2)\in \mathbb{N}^2$,
$y^\alpha=y_1^{\alpha_1} y_2^{\alpha_2}$, then:
\begin{gather*}
D^\alpha G_j(\omega_1,\omega_2)=\frac{1}{(2\pi)^2}\int\limits
e^{i(y_1\omega_1+y_2\omega_2)}y^\alpha \hat{G_j}(y_1,y_2)dy\\
=\frac{1}{(2\pi)^2}\int\limits
e^{i(y_1(\omega_1-x_1)+y_2(\omega_2-x_2))}y^\alpha G_j(x_1,x_2)dy
dx.
\end{gather*}
For this integral we can inverse our diffeomorphism via
$x_1(r,\theta)=r$ and $x_2(r,\theta)=p(\theta)$, locally on
$\mathrm{supp}(\Omega_j)$. For the $r$-integration we can extend
the integrand by 0 for $r<0$ and we obtain first:
\begin{gather*}
D^\alpha G_j(\omega_1,\omega_2)=\frac{1}{(2\pi)}\int\limits
e^{i\langle y_2,\omega_2-p(\theta)\rangle}y_2^{\alpha_2}
\frac{\partial}{\partial \omega_1 ^{\alpha_1} }\int
\Omega_j(\theta)
g(\omega_1\theta) \omega_1^{n-1}  d\theta dy_2,\\
=\frac{1}{(2\pi)}\int\limits e^{i\langle
y_2,\omega_2-p(\theta)\rangle}y_2^{\alpha_2}
\frac{\partial}{\partial \omega_1 ^{\alpha_1} }\int \omega_2
\Omega_j(\theta) \hat{f}(\omega_1\theta) \omega_1^{k+n-1}  d\theta
dy_2.
\end{gather*}
The remaining integral is simply the exterior derivative, of order
$\alpha_2$, of the Liouville measure on the surface
$p(\theta)=\omega_2$. For $\alpha_2=2$, observe that:
\begin{gather*}
\mathfrak{L}^{(2)}(p(\theta)\Omega_j(\theta) \hat{f}(\omega_1
\theta)(\omega_2)=\frac{\partial^2}{\partial\omega_2^2}\left(
\omega_2
\mathfrak{L}(\Omega_j(\theta) \hat{f}(\omega_1 \theta))(\omega_2)\right)\\
=\omega_2 \mathfrak{L}^{(2)}(\Omega_j(\theta) \hat{f}(\omega_1
\theta))(\omega_2)+2 \mathfrak{L}^{(1)}(\Omega_j(\theta)
\hat{f}(\omega_1 \theta))(\omega_2),
\end{gather*}
and that by construction the functions
$\mathfrak{L}(\Omega_j\hat{f})$ are smooth. Choosing
$\alpha_1=2k$, we have obtained:
\begin{gather}
\frac{\partial^{2k+2} G_j}{\partial \omega_1^{2k}
\partial \omega_2^2} (\omega_1,\omega_2)=\frac{\partial^{2k}}{\partial
\omega_1^{2k}}\left( \omega_1^{k+n-1}
\mathfrak{L}^{(2)}(p(\theta)\Omega_j(\theta) \hat{f}(\omega_1
\theta)(\omega_2) \right ) \notag\\
=\frac{\partial^{2k}}{\partial \omega_1^{2k}}\left(
\omega_1^{k+n-1} (\omega_2\mathfrak{L}^{(2)}(\Omega_j(\theta)
\hat{f}(\omega_1 \theta)(\omega_2)+2
\mathfrak{L}^{(1)}(\Omega_j(\theta) \hat{f}(\omega_1
\theta)(\omega_2) \right ). \label{value amplitude}
\end{gather}
By degree considerations w.r.t. $\omega_1$ we have respectively:
\begin{gather*}
\int\limits_{\mathbb{R}}\log(|\omega_2|)) \frac{\partial^{2k+1}
G_j}{\partial \omega_1^{2k-1}
\partial \omega_2^2} (0,\omega_2) d\omega_2=
\left\{
\begin{matrix}
 0\text{ if } k<n,\\
C(f)\neq 0 \text{ if } k\geq n.
\end{matrix}
\right.\\
\int\limits_{\mathbb{R}} \log(|\omega_2|)^2 \frac{\partial^{2k+1}
G_j}{\partial \omega_1^{2k-1}
\partial \omega_2^2} (0,\omega_2) d\omega_2=
\left\{
\begin{matrix}
 0\text{ if } k<n,\\
D(f)\neq 0 \text{ if } k\geq n.
\end{matrix}
\right.
\end{gather*}
Where $C$ and $D$ are obtained by inserting Eq.(\ref{value
amplitude}) in the integrals. Finally, in Eq.(\ref{value M''}) the
term attached to the product of logarithms is given by:
\begin{gather*}
\int\limits_{\mathbb{R}^{+}\times\mathbb{R}}\log(\omega_1)\log(|\omega_2|)
\frac{\partial^{2k+2} G_j}{\partial \omega_1^{2k}
\partial \omega_2^2} (\omega_1,\omega_2) d\omega_1 d\omega_2,\text{ if $k\geq n$},\\
(2k-1)! \int\limits_{\mathbb{R}^{+}\times
\mathbb{R}}\log(|\omega_2|)\frac{\partial^{2} G_j}{\partial
\omega_2^2} (\omega_1,\omega_2)
\frac{d\omega_1}{\omega_1^{2k}}d\omega_2, \text{ if } k<n.
\end{gather*}
For $k\geq n$ integrations by parts are not allowed but we can
anyhow conclude with Eq.(\ref{value amplitude}). We treat now
separately parts
\textbf{A)} and \textbf{B)} of Theorem \ref{main}.\medskip\\
\textit{Proof of part \textbf{A)}.}\\
To obtain the final result we sum over the partition of unity.
According to the considerations of homogeneity above, for $k<n$
the full contribution is generated by $\mu_0^{\mathrm{ell}}$ and
$M_j''(0)$. With the explicit values of $h(0)$ and $m(0)$, we
obtain that $(2\pi)^n\mu_0(f)$ equals:
\begin{equation*}
\int\limits_{\mathbb{R}^{+}\times\mathbb{S}^{n-1}}
\Omega_0(\theta) \hat{f}(r\theta) r^{n-k-1}
\frac{d\theta}{p(\theta)} dr +\sum\limits_{j}
\int\limits_{\mathbb{R}^+\times\mathbb{R}}
\log(|\omega_2|)\frac{\partial^{2} G_j}{\partial \omega_2^2}
(\omega_1,\omega_2)\frac{d\omega_1}{\omega_1^{2k}} d\omega_2.
\end{equation*}
With $\Omega_0=0$ near $\mathfrak{C}(p)\cap \mathbb{S}^{n-1}$, we
have $\mathfrak{L}(\Omega_0(\theta)\hat{f}(r\theta))(u)=0$ in a
neighborhood of $u=0$. Hence, in the first term, the integral
w.r.t. $\theta$ equals:
\begin{equation*}
\int\limits_{u\in\mathbb{R}}
\mathfrak{L}(\Omega_0(\theta)\hat{f}(r\theta))(u)
\frac{du}{u}=\left\langle \log(|u|)\, ;\,
\mathfrak{L}^{(1)}(\Omega_0(\theta)\hat{f}(r\theta))(u)
\right\rangle.
\end{equation*}
The derivation is in sense of distributions. For the coefficients
attached to $M_j''(0)$ we obtain:
\begin{gather*}
\int\limits_{\mathbb{R}}u\log(|u|)\mathfrak{L}^{(2)}(\Omega_j(\theta)\hat{f}(r\theta))(u)
 du + 2\int\limits_{\mathbb{R}}
\log(|u|)\mathfrak{L}^{(1)}(\Omega_j(\theta)\hat{f}(r\theta))(u)
du.
\end{gather*}
Since $(u\log(|u|)'=\log(|u|)+1$, via one integration by parts:
\begin{equation*}
\int\limits_{\mathbb{R}}u\log(|u|)\mathfrak{L}^{(2)}(\Omega_j(\theta)\hat{f}(r\theta))(u)
 du
=-\int\limits_{\mathbb{R}}\mathfrak{L}^{(1)}(\Omega_j(\theta)\hat{f}(r\theta))(u)
(\log(|u|)+1) du.
\end{equation*}
Observe the minus sign which fits with the weak derivation above.
Since for each $r$ and $j>0$,
$u\mapsto\mathfrak{L}(\Omega_j(\theta)\hat{f}(r\theta))(u)\in
C_0^{\infty}(\mathbb{R})$, we get:
\begin{equation*}
\int\limits_{\mathbb{R}}\mathfrak{L}^{(1)}(\Omega_j(\theta)\hat{f}(r\theta))(u)
du=0.
\end{equation*}
By integration w.r.t. $r$ and summation over the partition of
unity we obtain:
\begin{equation*}
\mu_0(f)= \frac{1}{(2\pi)^n}\int\limits_{0}^\infty \left\langle
\log(|\omega_2|)\, ; \, \mathfrak{L}^{(1)}
(\hat{f}(r\theta))(\omega_2) \right\rangle r^{n-k-1} dr,
\end{equation*}
which is the desired result when $k<n$.\medskip\\
\textit{Proof of part \textbf{B)}.}\\
Now, we consider $k\geq n$. All coefficients contribute via:
\begin{gather*}
(2\pi)^n\mu_0(f)=-h'(0)\frac{\partial ^{2k-1}}{\partial r^{2k-1}}
\left(r^{n+k-1} \int\limits_{\mathbb{S}^{n-1}} \hat{f}(r\theta)
\Omega_0(\theta)\frac{d\theta}{p(\theta)}
\right)_{|r=0}\\
+2kh(0)\int\limits_{\mathbb{R}^{+}\times\mathbb{S}^{n-1}}\log(r|p(\theta)|^{\frac{1}{k}})
\frac{\partial^{2k}} {\partial r^{2k}} \left(\hat{f}(r\theta)
r^{n+k-1}\right)
\Omega_0(\theta)\frac{d\theta}{p(\theta)}  dr\\
+ \frac{1}{2} \sum\limits_{j} (m(0)M''_j (0) +2m'(0)M_j'(0)).
\end{gather*}
If we split the integral with the logarithm we obtain two terms:
\begin{gather*}
\int\limits_{\mathbb{R}^{+}\times \mathbb{S}^{n-1}}\log(r)
\frac{\partial^{2k}} {\partial r^{2k}} (\hat{f}(r\theta)
r^{n+k-1})\Omega_0(\theta) \frac{d\theta}{p(\theta)}dr\\
-\frac{1}{k} \frac{\partial^{2k-1}} {\partial r^{2k-1}}\left(\,
\int\limits_{\mathbb{S}^{n-1}}\hat{f}(r\theta)
\log(|p(\theta)|)\Omega_0(\theta)\frac{d\theta}{p(\theta)}\right)_{|r=0}
.
\end{gather*}
Observe that, by construction, all integrals are well defined.
First, we express the contributions near $\mathfrak{C}(p)$.
Combining Eq.(\ref{value M'}) and Eq.(\ref{value amplitude}), we
find that:
\begin{gather*}
M_j'(0)=2\int\limits_{\mathbb{R}}\log(|\omega_2|))
\frac{\partial^{2k-1}}{\partial \omega_1^{2k-1}}\left(
\omega_1^{k+n-1} (\omega_2\mathfrak{L}^{(2)}(\Omega_j(\theta)
\hat{f}(\omega_1 \theta)(\omega_2)\right)_{|\omega_1=0} d\omega_2\\
+4\int\limits_{\mathbb{R}}\log(|\omega_2|))
\frac{\partial^{2k-1}}{\partial
\omega_1^{2k-1}}\left(\omega_1^{k+n-1}
\mathfrak{L}^{(1)}(\Omega_j(\theta) \hat{f}(\omega_1
\theta)(\omega_2) \right )_{|\omega_1=0}
 d\omega_2.
\end{gather*}
This term can be treated as in part \textbf{A)} and we obtain:
\begin{equation*}
M_j'(0)= 2\int\limits_{\mathbb{R}}\log(|\omega_2|))
\frac{\partial^{2k-1}}{\partial
\omega_1^{2k-1}}\left(\omega_1^{k+n-1}
\mathfrak{L}^{(1)}(\Omega_j(\theta) \hat{f}(\omega_1
\theta)(\omega_2) \right )_{|\omega_1=0}
 d\omega_2.
\end{equation*}
Next, combining Eq.(\ref{value M''}) and Eq.(\ref{value
amplitude}) we have:
\begin{gather*}
M_j''(0)=4 \int\limits_{\mathbb{R}} \log(|\omega_2|)^2
\frac{\partial^{2k-1}}{\partial \omega_1^{2k-1}}\left(
\omega_1^{k+n-1} (\omega_2\mathfrak{L}^{(2)}(\Omega_j(\theta)
\hat{f}(\omega_1 \theta)(\omega_2)\right)_{|\omega_1=0} d\omega_2\\
+8 \int\limits_{\mathbb{R}} \log(|\omega_2|)^2
\frac{\partial^{2k-1}}{\partial \omega_1^{2k-1}}\left(
\omega_1^{k+n-1} \mathfrak{L}^{(1)}(\Omega_j(\theta)
\hat{f}(\omega_1 \theta)(\omega_2)
\right )_{|\omega_1=0} d\omega_2\\
+8k\int\limits_{\mathbb{R}^{+}\times\mathbb{R}}\log(\omega_1)\log(|\omega_2|)
\frac{\partial^{2k}}{\partial \omega_1^{2k}}\left(
\omega_1^{k+n-1} (\omega_2\mathfrak{L}^{(2)}(\Omega_j(\theta)
\hat{f}(\omega_1 \theta)(\omega_2)\right) d\omega_1 d\omega_2\\
+16k
\int\limits_{\mathbb{R}^{+}\times\mathbb{R}}\log(\omega_1)\log(|\omega_2|)
\frac{\partial^{2k}}{\partial \omega_1^{2k}}\left(
\omega_1^{k+n-1} \mathfrak{L}^{(1)}(\Omega_j(\theta)
\hat{f}(\omega_1 \theta)(\omega_2) \right ) d\omega_1 d\omega_2.
\end{gather*}
The last two integrals can be combined as above. For the others,
we use:
\begin{equation*}
(u\log(|u|)^2)'= \log(|u|)^2+ 2\log(|u|),\, \forall u\neq 0,
\end{equation*}
and proceed to integrations by parts, which is legal since the
factors $\mathfrak{L}^{(k)}(.)(\omega_2)$ vanish for $\omega_2$
large and $\omega_2 \log(|\omega_2|)$ also vanishes at the origin.
We obtain:
\begin{gather*}
M_j''(0)=4 \int\limits_{\mathbb{R}} \log(|\omega_2|)^2
\frac{\partial^{2k-1}}{\partial \omega_1^{2k-1}}\left(
\omega_1^{k+n-1} \mathfrak{L}^{(1)}(\Omega_j(\theta)
\hat{f}(\omega_1 \theta)(\omega_2)
\right )_{|\omega_1=0} d\omega_2\\
-8\int\limits_{\mathbb{R}}\log(|\omega_2|)
\frac{\partial^{2k-1}}{\partial \omega_1^{2k-1}}\left(
\omega_1^{k+n-1} \mathfrak{L}^{(1)}(\Omega_j(\theta)
\hat{f}(\omega_1 \theta)(\omega_2)
\right )_{|\omega_1=0} d\omega_2\\
+8k\int\limits_{\mathbb{R}^{+}\times\mathbb{R}}\log(\omega_1)\log(|\omega_2|)
\frac{\partial^{2k}}{\partial \omega_1^{2k}}\left(
\omega_1^{k+n-1} (\mathfrak{L}^{(1)}(\Omega_j(\theta)
\hat{f}(\omega_1 \theta)(\omega_2)\right) d\omega_1 d\omega_2.
\end{gather*}
Observe that we have 3 different coefficients, like for the
coefficients attached to the set $\Omega_0$. We combine each of
these contributions by nature and by gathering carefully the
constants. First, we consider the term involving two logarithms:
\begin{gather*}
2kh(0)\int\limits_{\mathbb{R}^{+}\times \mathbb{S}^{n-1}}\log(r)
\frac{\partial^{2k}} {\partial r^{2k}} (\hat{f}(r\theta)
r^{n+k-1}) \Omega_0(\theta)\frac{d\theta}{p(\theta)}
dr\\+4km(0)\sum\limits_{j}\int\limits_{\mathbb{R}^{+}\times\mathbb{R}}\log(\omega_1)\log(|\omega_2|)
\frac{\partial^{2k}}{\partial \omega_1^{2k}}\left(
\omega_1^{k+n-1} \mathfrak{L}^{(1)}(\Omega_j(\theta)
\hat{f}(\omega_1
\theta))(\omega_2)\right) d\omega_1 d\omega_2\\
=\frac{1}{(k-1)!}\int\limits_{\mathbb{R}^{+}}
\log(\omega_1)\frac{\partial^{2k}}{\partial \omega_1^{2k}} \left(
\omega_1^{k+n-1} \left\langle \log(|\omega_2|)\, ; \,
\mathfrak{L}^{(1)}(\hat{f}(\omega_1
\theta))(\omega_2)\right\rangle \right) d\omega_1.
\end{gather*}
The change of sign for comes from a derivation in the sense of
distributions, a similar comment applies below. Next, we have:
\begin{gather*}
-\frac{2kh(0)}{k}\int\limits_{\mathbb{S}^{n-1}} \log(|p(\theta)|))
\frac{\partial^{2k-1}} {\partial r^{2k-1}} (\hat{f}(r\theta)
r^{n+k-1})_{|r=0} \Omega_0(\theta)\frac{d\theta}{p(\theta)}
\\
+2m(0)\sum\limits_{j}\int\limits_{\mathbb{R}}\log(|\omega_2|)^2
\frac{\partial^{2k-1}}{\partial \omega_1^{2k-1}}\left(
\omega_1^{k+n-1} \mathfrak{L}^{(1)}(\Omega_j(\theta)
\hat{f}(\omega_1
\theta))(\omega_2)\right)_{|\omega_1=0} d\omega_2\\
=\frac{1}{(2k)!}\frac{\partial^{2k-1}}{\partial \omega_1^{2k-1}}
\left( \omega_1^{k+n-1} \left\langle \log(|\omega_2|)^2 \, ; \,
\mathfrak{L}^{(1)}(\hat{f}(\omega_1
\theta))(\omega_2)\right\rangle \right)_{|\omega_1=0}.
\end{gather*}
Finally, we combine the remaining terms to obtain:
\begin{gather*}
-h'(0)\int\limits_{\mathbb{S}^{n-1}}\frac{\partial^{2k-1}}{\partial
r^{2k-1}} \left( r^{n+k-1} \hat{f}(r\theta)\right)_{|r=0}
\Omega_0(\theta)\frac{d\theta}{p(\theta)}
\\
+(2m'(0)-4m(0))\sum\limits_{j}\int\limits_{\mathbb{R}}\log(|\omega_2|)
\frac{\partial^{2k-1}}{\partial \omega_1^{2k-1}}\left(
\omega_1^{k+n-1} \mathfrak{L}^{(1)}(\Omega_j(\theta)
\hat{f}(\omega_1
\theta))(\omega_2)\right)_{|\omega_1=0} d\omega_2\\
=\frac{\gamma+\Psi(k)}{\Gamma(2k)}\frac{\partial^{2k-1}}{\partial
\omega_1^{2k-1}} \left( \omega_1^{k+n-1} \left\langle
\log(|\omega_2|)\, ; \, \mathfrak{L}^{(1)}(\hat{f}(\omega_1
\theta))(\omega_2)\right\rangle \right)_{|\omega_1=0}.
\end{gather*}
This proves parts \textbf{B)} of Theorem \ref{main}.\medskip\\
\textbf{Proof of Corollary \ref{coro main}}\medskip\\
We start by the analytic continuation of the elliptic part
$\mu^{\mathrm{ell}}(\zeta)$. The pole $\zeta=-1$ is simple and the
term of interest is given by:
\begin{equation*}
\mu_{0,-1}^{\mathrm{ell}}=\lim\limits_{\zeta\rightarrow 0} (\zeta
\mu^{\mathrm{ell}}(\zeta-1))
=\frac{1}{\Gamma(2k+1)}\partial_{\omega_1}^{2k-1} G(0).
\end{equation*}
The value of this coefficient was determined in the proof of
Theorem \ref{main}.

As concerns the singular term
$\mu^{\mathrm{sing}}(\zeta)=\sum\limits_{j=1}^{N}\mu_j^{\mathrm{sing}}(\zeta)$,
$\zeta=-1$ is a pole of order 2. Accordingly, the coefficients of
degree -2 and -1 are respectively given by:
\begin{gather*}
a_{-2,j}^{\mathrm{sing}}= \lim\limits_{\zeta\rightarrow 0} (\zeta^2\mu_j^{\mathrm{sing}}(\zeta-1)),\\
a_{-1,j}^{\mathrm{sing}}=\lim\limits_{\zeta\rightarrow 0}
\frac{\partial}{\partial \zeta}
(\zeta^2\mu_j^{\mathrm{sing}}(\zeta-1)).
\end{gather*}
Since $M_j(0)=0$, we have $a_{-2,j}^{\mathrm{sing}}= 0$ and
$a_{-1,j}^{\mathrm{sing}}= m(0)M'_j(0)$. To evaluate this
distributional coefficient we proceed exactly as above and obtain:
\begin{equation*}
a_{-1,j}^{\mathrm{sing}}= \frac{1}{\Gamma(1+2k)}
\int\limits_{\mathbb{R}}\log(|\omega_2|)) \frac{\partial^{2k+1}
G_j}{\partial \omega_1^{2k-1}
\partial \omega_2^2} (0,\omega_2) d\omega_2.
\end{equation*}
The discussion concerning the value of this term, established in
the proof of Theorem \ref{main}, gives the announced
result.$\hfill{\blacksquare}$\medskip\\
\textit{Duality brackets.}\\
Condition $(\mathcal{H})$ only insures that the Liouville measure
is smooth in a neighborhood of the origin. But the distributions
$\log(|y|)^\alpha$, $\alpha>0$, are smooth away from the origin.
With a smooth cut-off $\chi$, supported in a neighborhood of the
origin, we write $\langle \log(|y|)^\alpha \, ; \,
\mathfrak{L}^{(p)}(f)(y)\rangle$ as:
\begin{equation*}
\langle \log(|y|)^\alpha \, ;\, \chi(y)
\mathfrak{L}^{(p)}(f)(y)\rangle+\langle \log(|y|)^\alpha \, ; \,
(1-\chi(y))\mathfrak{L}^{(p)}(f)(y)\rangle.
\end{equation*}
Away from the origin, we can integrate by part the logarithmic
distribution. On the other side, we use that $y\mapsto
\mathfrak{L}(f)(y)$ is smooth on $\mathrm{supp}(\chi)$ if this
support is chosen small enough. This duality bracket is well
defined since both distribution have disjoint singular support.

Finally, this construction is independent from the cut-off $\chi$
if $\mathrm{supp}(\chi)$ is small enough with respect to the
covering of $\mathfrak{C}(p)$ introduced before. Conversely, for
any covering of $\mathfrak{C}(p)$ chosen such that
$|p(\theta)|\leq \varepsilon$ on each $\mathrm{supp}(\Omega_j)$,
$j\geq 1$, there exists a cut-off $\chi$ with the previous
properties. Hence the final
value is independent from the choice of our partition of unity on $\mathbb{S}^{n-1}$.\medskip\\
\textbf{Comments.}
\begin{itemize}
\item The relation between special functions, in particular
$\Gamma$ and hypergeometric, and fundamental solutions has
attracted much attention by the past. That's why we have greatly
detailed the
coefficients appearing in our setting.%
\item Residuum, and poles, of meromorphic distributions play also
an important r\^ole in asymptotic expansion of oscillatory and
fiber integrals. For example, the value of $m''(0)$ is exactly:
\begin{equation*}
\frac{ 12+k(6\gamma (2+k\gamma)+k\pi^2)+6k
(\Psi(2k)(2+2k\gamma + k \Psi(2k))-k \Psi^{(1)}(2k))}%
{3\Gamma(1+2k)},
\end{equation*}
where $\Psi^{(1)}(\zeta)=\partial_\zeta \Psi(\zeta)$ is the
polygamma-function of order 1. Such a coefficient is useful to
compute the second term of the asymptotic expansion of oscillatory
integrals with phase $p(\xi)$ or $p(\xi)^2$. See \cite{WON}
for this point.%
\item The determination of Liouville measures, and a fortiori of
their exterior differentials, is generally not possible. In the
case of homogeneous singularity, the determination of these
measures is sometimes possible in terms of generalized elliptic
integrals. See \cite{Wag} or \cite{Cam1} for different
examples.%
\item The condition that $k\in\mathbb{N}$ can be relaxed. We can
consider operators with a singularity at the origin providing that
their symbols are regular enough. If $\alpha>1$ is the degree, a
similar proof holds by using the integer part $k=[\alpha]+1$. All
constants are well defined as analytic functions of $\alpha$ and
one has to replace the radial derivations by the action of some
pseudo-differential operators with homogeneous symbol. If
$\alpha\leq 1$ the symbol $p$ is generally not $C^1$ and our
approach fails.
\end{itemize}


\begin{thebibliography}{00}
\bibitem{Bjo} J.-E. Bj\"ork, \textit{Rings of differential operators}. North-Holland
Math.Library, 21 (1979).
\bibitem{Cam} B. Camus, Fundamental solutions of homogeneous elliptic differential
operators, Bulletin des Sciences Math\'ematiques \textbf{130}
(2006), no. 3, 264-268.
\bibitem{Cam1} B. Camus, Asymptotic approximation of
degenerate fiber integrals, Journal of Mathematical Analysis and
Applications \textbf{320} (2006) no. 2, 30-44.
\bibitem{G-C} I.M. Guelfand and G.E. Chilov, \textit{Les distributions}.
Collection Universitaire de Math\'ematiques, VIII Dunod, Paris
(1962).
\bibitem{HOR1} L. H{\"o}rmander, \textit{The analysis of linear partial differential operators} 1,2,3,4,
Springer-Verlag (1985)
\bibitem{Wag} P. Wagner, Fundamental solutions of real homogeneous
cubic operators of principal type in three dimensions, Acta
Mathematica \textbf{182} (1999), 283-300.
\bibitem{Wag1} P. Wagner, On the fundamental solutions of a class of elliptic quartic
operators in dimension 3, J. Math. Pures Appl.(9) \textbf{81}
(2002), no. 11, 1191-1206.
\bibitem{WON} R. Wong, \textit{Asymptotic approximations of integrals},
Academic Press Inc. (1989).
\end{thebibliography}
\end{document}